\documentclass[12pt,twoside]{article} 
\usepackage{amsmath}  
\usepackage{amssymb}  
\usepackage{enumerate,epsfig}  
\newtheorem{lemma}{Lemma}[section]  
\newtheorem{theorem}{Theorem}[section]

\newtheorem{prop}{Proposition}[section]

\newtheorem{property}{Property}  
\newcommand{\R}{ {\mathbb R} }

\makeatletter  
\@addtoreset{equation}{section}  
\makeatother  
  
\newcommand{\cqfd}{{\unskip\kern 6pt\penalty 500  
\raise -2pt\hbox{\vrule\vbox to 6pt{\hrule width 6pt  
\vfill\hrule}\vrule}\par}}  
  
\newcommand{\ind}{{\mathbb I}}

\begin{document}  
  
\title{  DiPerna-Lions flow for  relativistic particles   
  in an electromagnetic field}  
\author{P.-E. JABIN  $^1$  and   
N. MASMOUDI$^{2}$   \\ $^1$   CSCAMM and Dpt of Mathematics,\\ 
 Univ. of 
  Maryland, College Park,  MD 20742-4015\\  
 email:  pjabin@cscamm.umd.edu  \\   
$^2$ Courant Institute \\  
251 Mercer St, New York, NY 10012 \\  
email: masmoudi@cims.nyu.edu   
\thanks{{\sf N. Masmoudi} is partially supported by an NSF Grant DMS-0703145} }  
    
\date{}

  
\maketitle  
\bigskip

\section*{Abstract}  
  
We show the existence and uniqueness of a DiPerna-Lions flow  for relativistic particles subject to   
a Lorentz force   in  an   
electromagnetic field. The electric and magnetic fields  solves the linear Maxwell system in the void but for singular initial conditions. As the corresponding force field is only in $L^2$, we have to perform a careful analysis of the cancellations over a trajectory.

\section{Introduction}  
%
An important open problem in the theory of renormalized solutions   
of the Boltzmann equation is to prove global existence for the   
 Boltzmann-Maxwell system. Indeed, even though global renormalized solutions    
are known to exist for the  Boltzmann equation \cite{DL89am}  
 and global weak solutions   
are known to exist for the Vlasov-Maxwell system \cite{DL89cpam},  there is   
no such result for the  Boltzmann-Maxwell system.  The main reason   
is that the    
two methods of proof   turn  out to be   incompatible. Indeed, on the   
one hand, the   
existence proof for the  Boltzmann equation  is based on the   
theory of renormalized solution and renormalization on the limit   
equation   requires a   
minimal regularity on the vector fields. More precisely on needs   
$W^{1,p}$ regularity, $p\geq 1$, see \cite{DL}, and \cite{Am} for the extension to   
the $BV$ case (we also refer to \cite{DeL} for a very nice presentation of the main results).  On the other hand, the existence proof for the   
  Vlasov-Maxwell system  uses a weak compactness argument that only   
uses averaging lemma and the solution constructed are weak solutions.   
It seems unclear how to put together the two methods since   
for the Boltzmann-Maxwell system the vector field for the transport   
equation in the phase space $(x,v)$ is only known to be   
in $L^\infty(0,T; L^2(\R^3))$ and hence we can not use   
the renormalization techniques for the limit.  
 
But of course we are here ignoring   
the extra structure coming from the fact that the electro-magnetic fields $(E,B)$ solve   
Maxwell equations with some right-hand side.   
  Let us also remark that for   
the  Boltzmann-Maxwell system with long range interaction  (without cut-off), existence   
of renormalized solutions can be proved since one has   
only to renormalize the equation on a regularized approximation    
  (before the passage to the limit)  and then take advantage of   
the strong convergence of $f^n(t,x,v)$  in all variables to pass to the limit   
and recover a renormalized solution  (see \cite{Diogo12}).    
    
\smallskip  
  
The goal of this paper is to give a first step in the   
understanding of   this problem, namely we study the case   
where $(E,B)$ solve the homogeneous  Maxwell system and take advantage of  
the different speeds of propagation between the slow particles and the fields which propagate at the speed of light. 
   
More precisely, we study the dynamics of relativistic particles in a   
given electromagnetic field. If one denotes $f(t,x,v)$ the  
distribution function in phase space, then $f$ solves the kinetic  
equation  
\begin{equation}\begin{split}  
&\partial_t f +\frac{v}{\sqrt{1+|v|^2}}  
\cdot\nabla_x f+ (E+v\times B)\cdot\nabla_v f=0,\quad x,v\in\R^3,\\  
&f(t=0,x,v)=f^0(x,v).  
\end{split}  
\label{kinetic0}  
\end{equation}  
for given $E$ and $B$.  
  
Proving well posedness for Eq. \eqref{kinetic0} is a completely open  
problem if one does not have any regularity for $E$ and $B$, which is  
the case that we wish to consider here. Therefore we assume that both  
fields are solution to Maxwell equations in the void. As an example, we will   
just consider the case where   
\begin{equation}  
E(t,x)=\partial_t \left(t\,\int_{S^2} E_0(x+t\,\omega)\,d\omega\right),\quad   
B(t,x)=\partial_t \left(t\,\int_{S^2} B_0(x+t\,\omega)\,d\omega\right).  
\label{wavestruct}  
\end{equation}  
Note that if $E_0$ and $B_0$ are only bounded then $E$ and $B$ have no extra regularity.  
 Nevertheless it is possible to show that  
\begin{theorem}  
Assume that $E_0,\;B_0\in L^1\cap L^{2}(\R^3)$ and $f^0\in  
L^1\cap L^\infty(\R^6)$.   
Then there exists a unique solution $f\in  
L^\infty(\R_+,\;L^p(\R^6))$ to Eq. \eqref{kinetic0}  
\label{ondes}  
\end{theorem}  
In our context, the first key point is that the particles solve  a kinetic equation (instead of a first order general transport equation). It has long been recognized that well posedness for kinetic equations is easier than for general transport equations. The $BV$ case was obtained earlier in \cite{Bo} (even improved to $BV_{loc}$ in \cite{Ha2}). While, in the general case, it has been shown that  the assumption of $BV$ regularity is optimal (see the counterexample in \cite{DeP}), it is possible to assume less for equations like \eqref{kinetic0}: $H^{3/4}$ for instance in \cite{CJ}.       
 
However in more than dimension $1$ ({\em i.e.} dimension 2 in phase space), $H^{3/4}$ is the best that can be done for the moment. Here it is hence necessary to use the additional structure provided on the fields by \eqref{wavestruct}. 
 
Of course the result relies on the fact that the particles move with  
a velocity which is strictly less than the speed of light and hence  
the wave equation possesses some regularizing properties. This idea   
was already used in the existence theory of strong solutions to the   
Vlasov-Maxwell system \cite{KS02,BGP03,GS86}, also this idea is at   
the origin of the space  resonance method used for instance in \cite{GMS11prep}.      
While this is perfectly satisfactory for the application we have in mind,   
one could nevertheless wish to study the interaction of particles with  
fields that are propagated at several different speeds, possibly  
comparable to the particles' velocity.  
 
We offer a partial answer in an  
essentially $1-d$ setting in~$x$ and now turn to  
\begin{equation}\begin{split}  
&\partial_t f +\alpha(v)  
\partial_x f+ F(t,x)\cdot\nabla_v f=0,\quad x\in \R, v\in \R^d,\\  
&f(t=0,x,v)=f^0(x,v).  
\end{split}  
\label{kinetic1}  
\end{equation}  
The function $\alpha$ is assumed to be Lipschitz: $\alpha\in  
W^{1,\infty}(\R^d)$ and satisfies a genuine non linearity assumption: There exists $C$ such that for all $w\in \R$ and   
$\eta > 0, $ we have   
\begin{equation} \label{alp-prop}  
|\{v, \ |\alpha(v) - w  |  \leq \eta \,    \} | \leq C \eta.   
\end{equation}  
  
The force field $F$ is assumed to be given by  
\begin{equation}  
F(t,x)=\sum_n F^0(x-\xi_n\,t)\,\mu_n,\label{structF}  
\end{equation}  
with $F^0 \in L^\infty$,   
together with the bound on the $\mu_n$  
\begin{equation}  
\exists \gamma>2,\quad \sum_n (1+n^\gamma)\,\mu_n<\infty.\label{mun}  
\end{equation}  
Note that whereas $x$ is necessarily $1$ dimensional, there is no such  
constraint on $v$. So for instance, Eq. \eqref{kinetic0} in the  
radially symmetric case would fit in this reduced framework.  
  
We have  
\begin{theorem} Assume \eqref{mun}, then for any $f^0\in L^p  
(\R\times\R^d)$ for $p>1$, there exists a unique solution $f\in  
  L^\infty(\R_+,\ L^p   
(\R\times\R^d)$ to Eq. \eqref{kinetic1} where $F$ is given by   
\eqref{structF}.\label{th1d}  
\end{theorem}  
In this 1 dimensional context, many results are already known. If $v\in \R$, $\alpha(v)=v$ and $F(t,x)=F^0(x)$ is autonomous then well posedness was already obtained in \cite{BD}, with an extension when $F^0$ is only $L^p$ in \cite{Ha1}. The key for both results is the Hamiltonian structure which implies the propagation of the total energy $v^2/2+\Phi(x)$ with $-\Phi'=F^0$ which allows to compute $v$ in terms of $x$ (up to a sign). 
 
This type of result was extended to general, autonomous transport equations in dimension 2 with a force field of bounded divergence, which is hence ``close enough'' in some sense to the Hamiltonian case. An additional assumption of noncharacteristic curve is also needed (this would correspond to \eqref{alp-prop} here); we refer to \cite{CR}, \cite{CCR} and \cite{ABC} which has the most general assumptions.     
 
Eq. \eqref{kinetic1} is still a kinetic equation: Even though strictly speaking, we are not in a Hamiltonian case as $v\in \R^d$, it is very close to the earlier formulation of \cite{BD} or \cite{Ha1} (more than the later extensions). But Theorem \ref{th1d} is not limited to autonomous $F$ which is the real improvement here. Unfortunately it is still not as general as one would like as we still have to assume some structure on the time dependence of $F$ given by \eqref{structF}-\eqref{mun}.  
\section{Proof of Theorem  \ref{ondes}}  
\subsection{Definition of the functional and reduction of the problem}  
As the structure of $E$ and $B$ is essentially the same, we only study  
the following equation  
\begin{equation}  
\partial_t f + \frac{v}{\sqrt{1+v^2}}\cdot\nabla_x  
f+F(t,x,v)\cdot\nabla_v f=0,\\  
\label{mainwave}  
\end{equation}  
with  
\begin{equation}  
F(t,x,v)=\nu(v)\,\partial_t \left(t\,\int_{S^2} F_0(x+t\omega)\,d\omega\right),  
\label{wavestruct2}  
\end{equation}   
with $\nu$ a $C^\infty$ function of $v$.  
  
Following the now classical connections between transport equations  
and ordinary differential equations, described in \cite{DL} or  
\cite{Am}, Theorem \ref{ondes} is implied by the existence and uniqueness  
of flows to the ODE  
\begin{equation}  
\begin{split}  
&\frac{d}{dt} X(t,x,v)=\frac{V(t,x,v)}{\sqrt{1+V^2(t,x,v)}},\quad   
\frac{d}{dt} V(t,x,v)=F(t,X),\\  
& X(0,x,v)=x,\qquad V(0,x,v)=v.  
\end{split}  
\label{mainode}  
\end{equation}  
As flows the solutions are required to satisfy   
\begin{property}\label{prop:Hamilt}  
For any  
  $t\in\R$ the application  
  \begin{equation}  
    (x,v)\in\R^3\times\R^3\mapsto(X(t,x,v),V(t,x,v))\in\R^3\times\R^3  
   \label{invertible}  
  \end{equation}  
  is globally invertible and has Jacobian $1$ at {\bf almost every}  
  $(x,v)\in\R^3\times\R^3$.  
 It also defines a semi-group   
\begin{equation}\begin{aligned}  
&  \forall s,t\in\R,\qquad & X(t+s,x,v)=X(s,X(t,x,v),V(t,x,v)), \\  
 & \mbox{and}\qquad & V(t+s,x,v)=V(s,X(t,x,v),V(t,x,v)).  
\end{aligned}\label{semigroup}\end{equation}  
\end{property}  
Theorem \ref{ondes} is then equivalent to  
\begin{theorem}  
Assume that $F$ satisfies Equation \eqref{wavestruct2} with initial data  
$F_0\in L^1\cap L^{2} (\R^3)$.  
Then there exists a unique solution to  
\eqref{mainode} satisfying Property~\ref{prop:Hamilt}.  
\label{ondes2}\end{theorem}  
  
Our strategy here is to derive explicit quantitative estimates on the trajectory. A functional was introduced in \cite{CD} for that (see also \cite{Ja} for an extension). We use here the modified functional introduced in \cite{CJ} specifically for kinetic equations: For any compact domain $\Omega\subset\R^6$  
we  
look at   
\begin{multline*}  
  Q_\delta(T)=\int_\Omega\log\left(1+\frac{1}{|\delta|^2} \left[\left(\sup_{0\leq  
        t\leq T}|X(t,x,v)-X^\delta(t,x,v)|^2  
      \right.\right.\right. \\  
      \left.\left.\left. +\int_0^T|V(t,x,v)-V^\delta(t,x,v)|^2 \:dt\right) \wedge 1  \right]  \right)  
      \:dx \:dv  
\end{multline*}  
where  $(X,V)$ is a solution to   
\eqref{mainode}   (or a regularized version) satisfying Property \ref{prop:Hamilt}  
and $(X^\delta,V^\delta)$ is   
 either a solution to a regularized version of \eqref{mainode}    
 or verifies 
\begin{equation}  
\begin{split}  
& \exists  
   (\delta_1,\delta_2)\in\R^6\quad\mbox{with}\ |(\delta_1,\delta_2) |\leq  
  \delta,\\  
&\qquad\qquad\qquad\qquad  
  (X^\delta,V^\delta)(t,x,v)=(X,V)(t,x+\delta_1,v+\delta_2).    
\end{split}\label{propdelta}  
\end{equation}  
%
In the following, we will frequently abuse the notation $\delta$ for in fact $|\delta|$. 
  
Then Theorem  \ref{ondes} is  
implied by  
\begin{prop}  
For any $\Omega$ compact, any $F_0\in L^1\cap L^{2}(\R^3)$,   
there exists a function $\psi$ depending   
on $T$, $\Omega$ and $F^0$,    
such that  for any  
$(X,V)$ solution to \eqref{mainode} with $F$ given by   
Equation \eqref{wavestruct}, satisfying Property  \ref{prop:Hamilt},  
and $(X^\delta,V^\delta)$ satisfying \eqref{propdelta},  
\[  
Q_\delta(T)\leq T\,\psi(-\log |\delta|),  
\]\label{propondes}  
with   
\[  
\frac{\psi(\xi)}{\xi}\longrightarrow 0,\quad\mbox{as}\ \xi\rightarrow  
\infty.  
\]  
\end{prop}  
The connection between  
Proposition  \ref{propondes} and Theorem  \ref{ondes2}  is simple and we refer   
to \cite{CD} or \cite{CJ} for a detailed explanation. Note however  
that it is not possible to obtain a direct  
Lipschitz estimate here and it is  
necessary to distinguish between $X-X^\delta$ and $V-V^\delta$ as in $Q_\delta$.   
%
\subsection{Proof of Proposition  \ref{propondes}: First steps}  
\subsubsection{Truncation of large velocities}  
%
First note that, by the usual estimates on solutions to wave  
equations,  
 since $F_0$ is in $L^1\cap L^2$ then so is $F(t,x)$ and  
$B$. Indeed  
\[  
F(t,x)=\int_{S^2} F_0(x+t\omega)\,d\omega+t\,\int_{S^2}  
\omega\cdot\nabla_x  
F_0(x+t\omega)\,d\omega.  
\]   
The bound is obviously true for the first term. As for the second,  
applying Fourier transform in $x$ yields  
\[  
{\cal F}t\,\int_{S^2}  
\omega\cdot\nabla_x  
F_0(.+t\omega)\,d\omega= \hat F_0(\xi)\,t\int_{S^2} e^{it\,\xi\cdot  
  \omega} \omega\cdot\xi\,d\omega=\hat F_0(\xi)\,M_t(\xi).  
\]  
The multiplier $M_t(\xi)$ is bounded uniformly in $t$.  
  
This means that for any $K$,   
\[\begin{split}   
&\int_{\Omega} \ind_{  \{ (x,v) |    \int_0^T |F(t,X(t,x,v))|\,dt\geq K   \}  }\,dx\,dv \\  
&\qquad\qquad\leq \frac{1}{K^2}\,  
\int_{\Omega}\int_0^T |F(t,X(t,x,v))|^2\,dt\,dx\,dv\leq C\,\frac{T}{K^2}.  
\end{split}  
\]  
Therefore denoting by  
$\Omega_K$   
\begin{equation}\begin{split}  
\Omega_K=\Bigg\{(x,v)\in \Omega\ s.t.&\quad \int_0^T |F(t,X(t,x,v))|\,dt\leq  
K\\  
& \mbox{and}\quad \int_0^T |F(t,X^\delta(t,x,v))|\,dt\leq  
K \Bigg\},  
\end{split}\end{equation}  
one deduces that $|\Omega\setminus\Omega_K|\leq CT/K^2$ and hence  
\begin{equation}  
Q_\delta(T)\leq \frac{CT}{K^2}\,|\log \delta|+Q_\delta^K(T),  
\end{equation}  
with  
\begin{multline*}  
Q_\delta^K(T)=\int_{\Omega_K}\log\left(1+\frac{1}{|\delta|^2}\left(\sup_{0\leq  
        t\leq T}|X(t,x,v)-X^\delta(t,x,v)|^2  
      \right.\right. \\  
      \left.\left. +\int_0^T|V(t,x,v)-V^\delta(t,x,v)|^2 \:dt\right)  \right)  
      \:dx \:dv.   
\end{multline*}  
Finally note that for $(x,v)$ in $\Omega_K$  
\[  
|V(t,x,v)|\leq |v|+K,\quad |X(t,x,v)|\leq  
|x|+t\,|v|+t\,K.   
\]  
As $(x,v)\in \Omega$ which is compact then for some constant $C$,  
$|V|+|X|\leq C\,K$ and the same is of course true for $X^\delta$ and  
$V^\delta$.   
%
\subsubsection{The free transport contribution}  
Now let  
\begin{equation*}\begin{split}  
  A_\delta(t,x,v)=&|\delta|^2+\sup_{0\leq  
    s\leq t}|X(s,x,v)-X^\delta(s,x,v)|^2\\  
&  +\int_0^t|V(s,x,v)-V^\delta(s,x,v)|^2 \:ds.  
\end{split}\end{equation*}  
Compute  
\begin{align*}  
  & \frac{d}{dt}\log\Biggl(1+\frac{1}{|\delta|^2}\Big(\sup_{0\leq s\leq  
        t}|X(s,x,v)-X^\delta(s,x,v)|^2\\  
&\qquad\qquad\qquad\qquad\qquad    
+\int_0^t|V(s,x,v)-V^\delta(s,x,v)|^2 \:ds\Big)\Bigg) \\  
  &    
  =\frac{2}{A_\delta(t,x,v)}\Bigg(\frac{d}{dt}\left(\sup_{0\leq s\leq  
      t}|X(s,x,v)-X^\delta(s,x,v)|^2\right)  
+(V(t)\!-\!V^\delta(t))\\  
&\qquad\times\int_0^t(F(s,(X,V)(s,x,v))  
  -F(s,(X^\delta,V^\delta)(s,x,v))) \:ds\Bigg).  
\end{align*}

Recall that for  any $f\in BV(0,T)$, we have   
\begin{equation*}  
  \frac{d}{dt}\left(\max_{0\leq s\leq t} f(s)^2\right)  
  \leq 2|f(t)f'(t)|\leq 4|f(t)|^2+\frac{1}{2}|f'(t)|^2.   
\end{equation*}  
And in addition  
\[  
|\partial_t X-\partial_t X^\delta|^2=|V/\sqrt{1+|V|^2}  
-V^\delta/\sqrt{1+|V^\delta|}|^2\leq 4|V-V^\delta|^2.  
\]  
Hence we deduce from the previous computation that  
\begin{align*}  
  Q_\delta^K(T) & \leq \int_{\Omega_K}\int_0^T  
  \frac{8\,|X-X^\delta|^2+|V-V^\delta|^2/2}  
  {A_\delta(t,x,v)} \:dt \:dx \:dv + \tilde{Q}_\delta(T) \\ & \leq  
  8|\Omega|(1+T)+\tilde{Q}_\delta(T)+\frac{1}{2}  
\int_{\Omega_K}\int_0^T  
  \frac{|V-V^\delta|^2}  
  {A_\delta(t,x,v)} \:dt \:dx \:dv  
\end{align*}  
where,  
\begin{multline*}  
  \tilde{Q}_\delta(T)=2\int_{\Omega_K}\int_0^T  
  \frac{V^\delta(t,x,v)-V(t,x,v)}{A_\delta(t,x,v)}\:\cdot \\  
  \int_{0}^{t}(F(s,(X^\delta,V^\delta)(s,x,v)-F(s,(X,V)(s,x,v)))  
\:ds \:dt\:dx\:dv.  
\end{multline*}  
Remark that  
\[\begin{split}  
\int_{\Omega_K}\int_0^T  
  \frac{|V-V^\delta|^2}  
  {A_\delta(t,x,v)} \:dt \:dx \:dv&\leq \int_{\Omega_K}\int_0^T  
  \frac{\partial_t A_\delta(t,x,v)}  
  {A_\delta(t,x,v)} \:dt \:dx \:dv\\  
&\leq \int_{\Omega_K} \log\left(\frac{A_\delta(T,x,v)}{|\delta|^2}\right)\:dx  
  \:dv\\  
&\leq Q_\delta^K(T),   
\end{split}\]  
where we recall that $ A_\delta(0,x,v) = \delta^2 $.   
Therefore, we have  
\begin{align*}  
  Q_\delta^K(T) & \leq   
  8|\Omega|(1+T)+2\,\tilde{Q}_\delta(T),  
\end{align*}  
and it is enough to bound $\tilde{Q}_\delta(T)$.   
   
For technical reasons related to some interpolations that will be explained later,  
 we will  bound a more general term than $\tilde  
Q_\delta$, namely     
\begin{multline*}  
  \bar{Q}_\delta(T)=2\int_{\Omega_K}\int_0^T  
  \frac{V\delta(t,x,v)-V(t,x,v)}{A_\delta(t,x,v)}\:\cdot \\  
  \int_{0}^{t}(\nu(V^\delta_s)\,G(s,X^\delta(s,x,v))-\nu(V_s)\,G(s,X(s,x,v)))  
\:ds \:dt\:dx\:dv,  
\end{multline*}  
where we assume that  $G$ satisfies  the same assumption  as  $F$, namely   
\begin{equation}  
G(t,x)=\partial_t\left(t\,\int_{S^{2}} G^0(x-\omega t)\,d\omega  
\right).\label{structG}   
\end{equation}   
In the term $\bar Q_\delta$ we decouple the connection between the dynamics  
of $(X,V)$ and  $(X^\delta, V^\delta)$ which is related to $F$   
  and  the   $G$ function which appears in $\bar Q_\delta$. This means that  
$\bar Q_\delta$ is now linear in $G^0$ which will later allow us to use  
interpolation theory.  
  
Let us remark  
that  
\[\begin{split}  
\bar{Q}_\delta(T)\leq &\int_{\Omega_K}\int_0^T   
(\nu(V)+\nu(V^\delta))\,  
\frac{V(t,x,v)-V^\delta(t,x,v)}{A_\delta(t,x,v)}\\  
&\qquad\quad\int_0^t  
(G(s,X^\delta_s)-G(s,X_s))\,ds\, dx\,dv\,dt\\  
&+\int_{\Omega_K}\int_0^T   
(\nu(V_t)-\nu(V^\delta_t))\,  
\frac{V(t,x,v)-V^\delta(t,x,v)}{A_\delta(t,x,v)}\\  
&\qquad\quad\int_0^t  
(G(s,X^\delta_s)+G(s,X_s))\,ds\, dx\,dv\,dt.  
\end{split}\]  
As $\nu$ is lipschitz the second term may be directly bounded by  
\[\begin{split}  
C_K\,\int_{\Omega_K}\int_0^T&   
\frac{|V(t,x,v)-V^\delta(t,x,v)|^2}{A_\delta(t,x,v)}\\  
&\quad\int_0^t  
(|G(s,X^\delta_s)|+|G(s,X_s)|)\,ds\, dx\,dv\,dt,  
\end{split}\]  
where $C_K=\max_{B(0,C\,K)} |\nabla \nu(K)|$. Note that  
\[  
\int_s^T \frac{|V(t,x,v)-V^\delta(t,x,v)|^2}{A_\delta(t,x,v)}\,dt\leq  
\int_s^T \frac{\partial_t A_\delta}{A_\delta}\,dt\leq -C\,\log |\delta|.   
\]  
Hence by Fubini   
\begin{equation}  
\bar Q_\delta(T) \leq   
2\,I_\delta(T)+C_K\,T\,|\log  
  |\delta||\,(\|G^0\|_{L^1}+\|G^0\|_{L^2}),  
\end{equation}  
with  
\[\begin{split}  
I_\delta=\nu(K)\,&  
\int_{\Omega_K}\int_0^T   
\frac{|V(t,x,v)-V^\delta(t,x,v)|}{A_\delta(t,x,v)}\\  
&\quad\left|\int_0^t\,  
(G(s,X^\delta_s)-G(s,X_s))\,ds\,\right| dx\,dv\,dt.  
\end{split}\]  
%
\subsection{Proof of Proposition  \ref{propondes}: The main bound}  
The term $I_\delta$ is quite technical to bound and we hence summarize  
the computations in  
the following lemma  
\begin{lemma} For any $G$ satisfying \eqref{structG} with $G_0\in  
  L^1\cap L^{2}(\R^3)$   
and any $(X,V)$ solution to \eqref{mainode} with  
  \eqref{prop:Hamilt}, any $(X^\delta,V^\delta)$ verifying  
  \eqref{propdelta}, there exists a constant $C$  
depending  
  only  on $T$, s.t.  
\[  
I_\delta\leq C\,\nu(K)\,K^{10}\,|\log |\delta||\,(\|G^0\|_{L^1}  
+\|G^0\|_{L^2})\;(1+\|F^0\|_{L^2}).  
\]\label{teclem}    
\end{lemma}

{\bf Beginning of the proof of Lemma \ref{teclem}.} Write  
\[\begin{split}  
\int_0^t  
\big(G(s,X^\delta_s)-G(s,X_s)\big)\,ds&=\int_0^t\int_{S^{2}}  
\big(  G^0(X_s-\omega s)-G^0(X_s^\delta -\omega s) \big)\, d\omega ds\\  
+\int_0^t\int_{S^{2}}  
(\omega&\cdot\nabla_x  
 G^0(X_s-\omega  
s)-\omega\cdot\nabla_x G^0(X_s^\delta -\omega s))\,s\,d\omega ds.\\  
\end{split}\]  
Introduce the two changes of variables  
\[  
\Phi_X(s,\omega)=X_s-s\omega,\quad  
\Phi_{X^\delta}(s,\omega)=X_s^\delta-s\omega.  
\]  
The Jacobian of the transform is  
\[  
J_X=C\,s^2\,|\dot X_s\cdot \omega-1|,  
\]  
and the corresponding formula for $J_{X^\delta}$.   
Denote $(s_X,\omega_X)(z)$ the inverse of $\Phi_X$, namely  
$z = X_{s_X(z)} -s_X(z)  \omega_X(z) $  and  
$O_X^t=\bigcup_{s\in [0,\ t],\omega\in S^{2}} \Phi_X(s,\omega)$.   
One can easily prove that  
\begin{equation}  
O_X^t=B(X(t,x,v),t).\label{ball}  
\end{equation}  
One inclusion is indeed obvious and as for the other one, note that  
\[  
|\Phi_X(s,\omega)-X(t)|\leq s+|X(s)-X(t)|\leq s+|t-s|,  
\]  
as $|\dot X|<1$.  
  
One obtains  
\[  
| \int_0^t  
(G(s,X^\delta_s)-G(s,X_s))\,ds | \leq A+|B|+|C|+D+|E|,  
\]  
with  
\begin{equation}  
\begin{split}  
A&=\int_{O_X^t\setminus O_{X^\delta}^t}  
|G^0(z)| \;\frac{C\,dz}{s^2_{X}  
  |\dot X_s\cdot\omega_{X}-1|}\\  
&\qquad +\int_{O_{X^\delta}^t\setminus O_{X}^t}  
|G^0(z)| \;\frac{C\,dz}{s^2_{X^\delta}  
  |\dot X_s^\delta\cdot\omega_{X^\delta}-1|}\\  
B&=\int_{O_X^t\cap O_{X^\delta}^t} G^0(z)\,\left(\frac{C}{s_X^2  
  |\dot X_{s_X}\cdot\omega_X-1|}-\frac{C}{s_{X^\delta}^2  
  |\dot X_{s_{X^\delta}}^\delta\cdot\omega_{X^\delta}-1|}\right)\,dz,  
\end{split}\label{int}  
\end{equation}  
which correspond to the term without derivative,   
\begin{equation}  
\begin{split}  
C&=\int_{\partial B(X(t),t)} G^0(z)\,\frac{C}{s_X  
  |\dot X_{s_X}\cdot\omega_X-1|}\,dS(z)\\  
&\qquad-\int_{\partial B(X^\delta(t),t)}  
G^0(z)\, \frac{C}{s_{X^\delta}   
  |\dot X^\delta_{s_{X^\delta}}\cdot\omega_{X^\delta}-1|}\,dS(z),\\  
D&=\int_{O_X^t\setminus O_{X^\delta}^t}  
|G^0(z)|\, \left|  \nabla_z\cdot\left(\;\omega_X\,\frac{C}{s_{X}  
  |\dot X_s\cdot\omega_{X}-1|}\right) \right| \,dz\\  
&\qquad +\int_{O_{X^\delta}^t\setminus O_{X}^t}  
|G^0(z)|\; \left| \nabla_z\cdot\left(\omega_{X^\delta}\,\frac{C}{s^2_{X^\delta}  
  |\dot X_s^\delta\cdot\omega_{X^\delta}-1|}\right) \right| \,dz,\\  
\end{split}\label{int2}  
\end{equation}  
and finally  
\begin{equation}  
\begin{split}  
E&=\int_{O_X^t\cap O_{X^\delta}^t} G^0(z)\,\Bigg(\nabla_z\cdot\left(  
\frac{C\,\omega_X}{s_X  
  |\dot X_{s_X}\cdot\omega_X-1|}\right)\\  
&\qquad\qquad\qquad\qquad-\nabla_z\cdot\left(  
\frac{C\,\omega_{X^\delta}}{s_{X^\delta}  
  |\dot  
  X_{s_{X^\delta}}^\delta\cdot\omega_{X^\delta}-1|}\right)\Bigg)\,dz.\\  
\end{split}\label{int2'}  
\end{equation}  
Note that $C$ is a boundary term coming from the integration by parts.   
We hope there is no confusion since $C$ is used for constants that   
may change from one line to the other.   
We denote by $I_A$,\dots,$I_E$ the integrals, over  
$\Omega_K\times[0,\ T]$ of the previous  
quantities multiplied by  
\[  
(\nu(V_t)+\nu(V_t^\delta))\,\frac{|V(t,x,v)-V^\delta(t,x,v)|}{A_\delta(t,x,v)}.  
\]  
%
%
\subsubsection{Bound for $I_B$}  
%
As part of the bound for $I_E$ uses the same steps, we prove here a  
more general result on quantities like $I_B$. We call them  $I_{Bmod}$.  
\begin{lemma} Assume that $H\in W^{1,\infty}$ and that $4/3<p<2$,   
  then one has for some constant $C$ depending on $p$ and the norm of $H$   
and for $k=1,2$,   
\[\begin{split}  
I_{Bmod}:=&\nu(K)\,\int_0^T\int_{\Omega_K}  
\frac{|V_t-V^\delta_t|}  
{A_\delta(t,x,v)}   
\int_{O_X^t\cap O_{X^\delta}^t}  
|G^0(z)|\,\Bigg|\frac{H(s_X,\omega_X)}  
{s_X^2\,|\dot X_{s_X}^\delta\cdot\omega_X-1|^k}\\  
&\qquad\qquad -\frac{H(s_{X^\delta},\omega_{X^\delta})}  
{s_{X^\delta}^2\,|\dot X_{s_{X^\delta}}\cdot\omega_{X^\delta}-1|^k}  
\Bigg|\,dz\,dx\,dv\,dt\\  
&\leq C\,K^{5+k}\,\nu_K\,\sqrt{-\log  
  |\delta|}\,(\|G^0\|_{L^1}+\|G^0\|_{L^p}).    
\end{split}\]  
\label{lemIB}\end{lemma}  
  
Lemma \ref{lemIB} with  $k=1$ and  $H$ constant  implies that for $4/3<p\leq 2$, we have   
\begin{equation}  
I_B\leq C\,\nu(K)\,K^6\,\sqrt{-\log  
  |\delta|}\,(\|G^0\|_{L^1}+\|G^0\|_{L^p}).\label{IB0}   
\end{equation}  
  
\noindent{\bf Proof of Lemma \ref{lemIB}.} Denote  
\[  
Bmod=\int_{O_X^t\cap O_{X^\delta}^t}  
|G^0(z)|\,\Bigg|\frac{H(s_X,\omega_X)}  
{s_X^2\,|\dot X_{s_X}\cdot\omega_X-1|^k}
 -\frac{H(s_{X^\delta},\omega_{X^\delta})}  
{s_{X^\delta}^2\,|\dot X_{s_{X^\delta}}^\delta\cdot\omega_{X^\delta}-1|^k}  
\Bigg|\,dz.  
\]  
  
Recall that $\dot X_s=V_s/\sqrt{1+V_s^2}$ and that $|V_s|\leq C\,K$, so that  
$|\dot X_s|\leq /1-1/(CK)$ and $|\dot X_s\cdot \omega-1|\geq 1/(CK)$. Then  
\[\begin{split}  
\left|\frac{H(s_{X},\omega_{X})}{s_X^2  
  |\dot X_s\cdot\omega-1|^k}-\frac{H(s_{X^\delta},\omega_{X^\delta})}{s_{X^\delta}^2  
  |\dot X_s^\delta\cdot\omega_{X^\delta}-1|^k}\right|\leq  
&C\,K^k\,\Bigg(|s_X-s_{X^\delta}|\,\left(\frac{1}{s_X^3}  
+\frac{1}{s_{X^\delta}^3}\right) \\  
+K\,|\omega_X-\omega_{X^\delta}|\; \left(\frac{1}{s_X^2}  
+\frac{1}{s_{X^\delta}^2}\right)+K&\frac{|V_{s_X}-V_{s_X}^\delta|}{s_X^2}  
+K\,\frac{|V_{s_{X^\delta}}-V_{s_{X^\delta}}^\delta|}{s_{X^\delta}^2}\Bigg). \\  
\end{split}\]  
On the other hand, by definition  
$X_{s_X}-s_X\,\omega_X=X_{s_{X^\delta}}^\delta-s_{X^\delta}\,\omega_{X^\delta}$, so  
\begin{equation}  
\Phi_X(s_X,\omega_X)-\Phi_X(s_{X^\delta},\omega_{X^\delta})  
=X_{s_X}-X_{s_{X^\delta}}   
-s_X\,\omega_X+s_{X^\delta}\,\omega_{X^\delta}=  
X_{s_{X^\delta}}^\delta-   
X_{s_{X^\delta}}.\label{inter}  
\end{equation}  
Recall that $|\omega_X|=|\omega_{X^\delta}|=1$. Therefore  
$|s_X\,\omega_X-s_{X^\delta}\,\omega_{X^\delta}|\geq  
|s_X-s_{X^\delta}|$ and as $| X_{s_X}-X_{s_{X^\delta}} |\leq  
(1-1/(CK))|s_X-s_{X^\delta}|$,  
\[  
|s_X-s_{X^\delta}|\leq \max_{s\leq inf(s_X,s_{X^\delta})}  
|X_s-X_s^\delta|+(1-1/(CK))|s_X-s_{X^\delta}|,  
\]  
which implies  
\[  
|s_X-s_{X^\delta}|\leq CK\,\max_{s\leq inf(s_X,s_{X^\delta})}  
|X_s-X_s^\delta|.  
\]  
Using this estimate in \eqref{inter}, one concludes that  
\begin{equation}\begin{split}  
&|s_X-s_{X^\delta}|\leq C\,K\,\max_{s\leq inf(s_X,s_{X^\delta})}  |X_s-X_s^\delta|,\\  
&  
|\omega_X-\omega_{X^\delta}|\leq C\,K\,\max_{s\leq  
  inf(s_X,s_{X^\delta})}  
 |X_s-X_s^\delta|\;\left(\frac{1}{s_X}  
+\frac{1}{s_{X^\delta}}\right).  
\end{split}\label{deltasom}  
\end{equation}  
Inserting this in the term $I_{Bmod}$   
enables to bound it by  
\[\begin{split}  
I_{Bmod}\leq  
&C\,K^{k+2}\,\nu_K\int_0^T\int_{\Omega_K}  
\frac{|V(t,x,v)-V^\delta(t,x,v)|}{A_\delta(t,x,v)}\\  
& \Bigg(  
\int_{O_X^t\cap O_{X^\delta}^t} |G^0(z)|\max_{s\leq \inf(s_X,s_{X^\delta})}  
|X_s-X_s^\delta|\,\left(\frac{1}{s_X^3}  
+\frac{1}{s_{X^\delta}^3}\right)\;dz\\  
&+\int_{O_X^t\cap O_{X^\delta}^t}\,|G^0(z)|\left(  
\frac{|V_{s_X}-V_{s_X}^\delta|}{s_X^2}  
+\frac{|V_{s_{X^\delta}}-V_{s_{X^\delta}}^\delta|}{s_{X^\delta}^2}\right)\,dz\Bigg)  
\,dx\,dv\,dt,  
\end{split}\]  
with  
$\nu_K=\max_{B(0,K)} |\nu(v)|$.  
  
Changing back variables to $s,\omega$, one eventually finds   
\begin{equation}\begin{split}  
I_{Bmod}&\leq C\,K^{k+2}\,\nu_K\,  
\int_0^T\int_{\Omega_K}\frac{|V(t,x,v)-V^\delta(t,x,v)|}{A_\delta(t,x,v)}\\  
&\Bigg(\int_0^t\frac{\max_{r\leq s}   
|X_s-X_s^\delta|}{s}\,\int_{S^2}  
|G^0(X_s+s\omega)|\,d\omega\,ds\\  
&+  
\int_0^t\int_{S^2}  
|G^0(X_s+s\omega)|\,|V_s-V_s^\delta|\,d\omega\,ds+symmetric\ terms  
\Bigg)\,dx\,dv\,dt.  
\end{split}\label{boundB}\end{equation}  
Note that  
\[  
\max_{r\leq s}   
|X_s-X_s^\delta|\leq \sqrt{s}\,\left(\int_0^s  
|V_r-V_r^\delta|^2\,dr\right)^{1/2}\leq \sqrt{s}\,\sqrt{A_\delta(t,x,v)}.  
\]  
Hence by the definition of $A_\delta$  
\[\begin{split}  
&\int_0^T\int_{\Omega_K}\frac{|V-V^\delta|(t,x,v)}{A_\delta(t,x,v)}  
\int_0^t\frac{\max_{r\leq s}   
|X_s-X_s^\delta|}{s}\,\int_{S^2}  
|G^0(X_s+s\omega)|\,d\omega\,ds\\  
&\leq \int_0^T s^{-1/2}\,\int_{\Omega_K} \int_{S^2}  
|G^0(X_s+s\omega)|\,d\omega   
\int_s^T \frac{|V-V^\delta|(t,x,v)}{\sqrt{A_\delta(t,x,v)}}\,dt\,dx\,dv\,ds.   
\end{split}\]  
However as $|V-V^\delta|^2\leq \partial_t A_\delta$  
\begin{equation}\begin{split}  
\int_s^T \frac{|V-V^\delta|(t,x,v)}{\sqrt{A_\delta(t,x,v)}}\,dt&\leq  
\sqrt{T}\,\left(\int_s^T  
\frac{|V-V^\delta|^2(t,x,v)}{A_\delta(t,x,v)}\,dt\right)^{1/2}\\  
&\leq  
C\,\sqrt{-\log |\delta||},    
\end{split}\label{usefulest}\end{equation}  
one deduces that  
\[\begin{split}  
&\int_0^T\int_{\Omega_K}\frac{|V-V^\delta|(t,x,v)}{A_\delta(t,x,v)}  
\int_0^t\frac{\max_{r\leq s}   
|X_s-X_s^\delta|}{s}\,\int_{S^2}  
|G^0(X_s+s\omega)|\,d\omega\,ds\\  
&\leq C\, \sqrt{-\log |\delta||}\,\int_0^T s^{-1/2}\,\int_{\Omega_K} \int_{S^2}  
|G^0(X_s+s\omega)|\,d\omega\,dx\,dv\\  
&\leq C\, K^3\,\sqrt{-\log |\delta|}\,\|G^0\|_{L^1},  
\end{split}\]  
by a change of variables.  
  
Let us turn to the second term in\eqref{boundB}. Denote  
\[  
\tilde G(s,x)=\int_{S^2} |G^0(x+s\omega)|\,d\omega.  
\]  
By Cauchy-Schwartz, and since $A_\delta(t)\geq \int_0^t |V-V^\delta|^2$  
\[  
\begin{split}  
&\int_0^T\int_{\Omega_K}\frac{|V-V^\delta|(t,x,v)}{A_\delta(t,x,v)}  
\int_0^t\int_{S^2}  
|G^0(X_s+s\omega)|\,|V_s-V_s^\delta|\,d\omega\,ds\,dx\,dv\,dt\\  
&\quad\leq C\,\int_{\Omega_K}\left(\int_0^T |\tilde  
G(s,X_s)|^2\,ds\right)^{1/2}\,  
\int_0^T \frac{|V-V^\delta|(t,x,v)\,dt}{\sqrt{A_\delta(t,x,v)}}\,dx\,dv.  
\end{split}  
\]  
Hence  
\[  
\begin{split}  
&\int_0^T\int_{\Omega_K}\frac{|V-V^\delta|(t,x,v)}{A_\delta(t,x,v)}  
\int_0^t\int_{S^2}  
|G^0(X_s+s\omega)|\,|V_s-V_s^\delta|\,d\omega\,ds\,dx\,dv\,dt\\  
&\quad \leq C\,\sqrt{-\log|\delta|}\,\int_{\Omega_K}\left(\int_0^T |\tilde  
G(s,X_s)|^2\,ds\right)^{1/2}\,dx\,dv\\  
&\quad \leq C\,|\Omega|\,\sqrt{-\log|\delta|}\,K^3\,  
\int_0^T \int_{\R^3} |\tilde  
G(s,X_s)|^2\,dx\,ds.   
\end{split}\]  
By the usual estimates on solutions to wave equations (see for  
instance \cite{St}, chapter. 8, 5.21), one has  
\begin{lemma} For any $4/3  
<p<2$,  there exists   
$ C<\infty$ such that  for all $ G^0 \in L^p \cap L^1 $, we have   
\[  
\left\|\int_{S^2} G^0(x+s\omega)\,d\omega\right\|_{L^2}\leq  
C\,s^{3/2-3/p}\left(\|G^0\|_{L^1}+\|G^0\|_{L^p}\right).   
\]\label{wavedisp}  
\end{lemma}  
Combining all estimates and using Lemma \ref{wavedisp} (note that  
$3-6/p>-1$ if $p>4/3$), one finally obtains  
\begin{equation}  
I_{Bmod}\leq  
C\,K^{5+k}\,\nu_K\,\sqrt{-\log{|\delta|}}\,\left(\|G^0\|_{L^1}+\|G^0\|_{L^p}   
\right).\label{IB}  
\end{equation}  
\subsubsection{Bound on $I_A$}  
%
Note that $O^t_X\setminus  
  O^t_{X^\delta}=B(X_t,t)\setminus B(X_t^\delta,t)$.   
Since $|\dot X_s|\leq 1-1/CK$ with the same for $X^\delta$, one has that  
for any $\omega\in S^2$  
\[  
|X_t^\delta-X_s-s\omega|\leq |X_t^\delta-x|+|X_s-x|+s\leq t\,(1-1/CK)+2s<t,  
\]  
if $s<t/(2CK)$. Therefore   
\[  
\forall z\in O^t_X\setminus  
  O^t_{X^\delta},\quad s_X(z)\geq \frac{t}{CK},   
\]  
\[  
A\leq \frac{C\,K^2}{t^2}\,\int_{O^t_X\setminus  
  O^t_{X^\delta}}|G^0(z)|\,dz+ \frac{C\,K^2}{t^2}\int_{O^t_X\setminus  
  O^t_{X^\delta}}|G^0(z)|\,dz.   
\]  
Now we introduce the following modified maximal operator  
\[  
\tilde Mg(x)=\sup_{\eta\leq \varepsilon} \frac{1}{\varepsilon^2\,\eta}\,  
\int_{\varepsilon-\eta\leq |z-x|\leq \varepsilon} |g(z)|\,dz.  
\]  
Note that for example if $z\in O^t_{X}\setminus O^t_{X^\delta}$ then  
$|z-X_t|\leq t$  and $ t \leq |z-X_t^\delta|$. Using that   
 $|z-X_t^\delta|\leq |z-X_t | + |X_t-X_t^\delta|$, we deduce that  
\[  
t-|X_t-X_t^\delta|\leq |z-X_t|\leq t.  
\]  
Hence, one controls $A$ with $\tilde MG^0$ as  
\[  
A\leq C\,K^2\,|X_t-X_t^\delta|\,(\tilde MG^0(X_t)+\tilde MG^0(X_t^\delta)).  
\]  
This allows for an easy bound on $I_A$ in terms of the $L^1$ norm of  
$\tilde MG^0$  
\[\begin{split}  
I_A&=\int_0^T\int_{\Omega_K} (\nu(V_t)+\nu(V_t^\delta)\,  
\frac{|V_t-V_t^\delta|}{A_\delta(t,x,v)}\,A\,dx\,dv\,dt\\  
&\leq C\,K^2\,\nu(K)\,\int_0^T\int_{\Omega_K}  
\frac{|V_t-V_t^\delta|}{\sqrt{A_\delta(t,x,v)}}\,  
(\tilde MG^0(X_t)+\tilde MG^0(X_t^\delta))\,dx\,dv\,dt\\  
&\leq C\,K^5\,\nu(K)\,\sqrt{|\log \delta|}\, \left( \int_{|x|\leq C\,K}   
|\tilde MG^0(x)|^2 \right)^{1/2} ,   
\end{split}  
\]  
by the change of variable $(x,v)\rightarrow (X_t,V_t)$ (or  
$(X_t^\delta,V_t^\delta)$) and \eqref{usefulest}.  
  
We hence need to estimate $\tilde MG^0$. As it is defined, it turns out that  
it  is of the same order as the spherical maximal operator for which  
bounds are well known  
\begin{lemma} $\forall p>3/2$, $\exists C>0$ s.t. for any smooth function $g$  
\[  
\|\tilde Mg\|_{L^p(\R^3)}\leq C\, \|g\|_{L^p}.  
\]\label{modmaximal}  
\end{lemma}  
\noindent {\bf Proof of Lemma \ref{modmaximal}.} Simply note that  
\[  
\int_{\varepsilon-\eta\leq |z-x|\leq \varepsilon} |g(z)|\,dz\leq  
\int_{\varepsilon-\eta}^\varepsilon \int_{S(x,r)}  
|g(z)|\,dS(z)\,dr \leq C\,\eta\,\varepsilon^2\, M_S\,g(x),   
\]  
where $M_S$ is defined by  
\[  
M_S g(x)=\sup_r \oint_{S(x,r)}  
|g(z)|\,dS(z).  
\]  
It is proved (see \cite{St}) that in dimension $3$, $M_S$ is bounded on  
$L^p$ for any $p>3/2$ (the limit exponent is sharp) which easily implies the  
lemma.  
  
Note that this exponent is also sharp for $\tilde M$ as obviously  
\[  
\sup_{\eta\leq \varepsilon} \frac{1}{\eta}\int_{\varepsilon-\eta\leq  
  |z-x|\leq \varepsilon} |g(z)|\,dz\geq C\,\int_{S(x,\varepsilon)}  
|g(z)|\,dS(z).   
\]  
And so $\tilde Mg\geq C\, M_S g$.  
\cqfd  
  
Coming back to $I_A$   
one concludes that   
\begin{equation}  
I_A\leq C\,K^{8-3/2}\,  
\nu(K)\,\sqrt{|\log \delta|}\,\|G^0\|_{L^2}.\label{boundA}  
\end{equation}   
\subsubsection{Bound for $I_C$}  
Note that on $\partial B(X_t,t)$, one obviously has $s_X=t$ and  
$\omega_X=(z-X_t)/t$. Hence  
\[  
\int_{\partial B(X_t,t)} G^0(z)\,\frac{dS(z)}{s_X\,|\dot  
  X_{s_X}\cdot\omega_X-1|}  
= C\,t\,  
\int_{S^2} \frac{G^0(X_t+t\,\omega)}{|\dot X_t\cdot\omega-1|}\,d\omega.  
\]  
One finds that   
\[\begin{split}  
I_C&\leq C\,K^2\,\nu(K)\int_{\Omega_K}\int_0^T  
\frac{|V_t-V_t^\delta|^2}{A_\delta(t)}\,t\int_{S^2} |G^0(X_t+t\omega)|\,d\omega  
\,dx\,dv\,dt\\  
+& C\,\nu(K)\int_{\Omega_K}\!\!\int_0^T  
\frac{|V_t-V_t^\delta|}{A_\delta(t)}\,t\left|\int_{S^2}\frac{G^0(X_t+t\,\omega)  
-G^0(X_t^\delta+t\omega)}{|\dot X_t^\delta\cdot\omega-1|}\,d\omega\right|  
\,dx\,dv\,dt,\\  
&=I_C^1+I_C^2\end{split}\]  
Let us deal first with $I_C^1$. Noticing as before that  
\[  
\frac{|V_t-V_t^\delta|^2}{A_\delta(t)} \leq \partial_t \log A_\delta(t).  
\]  
Then denoting by $W$ the usual wave operator  
\[   
Wg (t,x)  =t\,\int_{S^2} g(x+t\omega)\,d\omega,  
\]   
one finds,  after integration by parts in time,   that  
\[\begin{split}  
I_C^1&\leq -C\,K^2\,\nu(K)\,\int_0^T\int_{\Omega_K} \log A_\delta \;\partial_t  
[W|G^0|(t,X_t)]\\  
&\quad+C\,K^2 \int_{\Omega_K} \log A_\delta(T,x,v)\,W|G^0|(T,X_T)\\  
&\leq C\,K^{13/2}\,\nu(K)\,|\log |\delta||\,(\|\partial_t  
W|G^0|\|_{L^2}+\|W|G^0|\|_{L^2}).  
\end{split}\]  
Of course $(\partial_t,\nabla_x) W|G^0|$   
is bounded  
in $L^2$ by the norm of $G^0$ in $L^2$ as can be seen by taking  
the Fourier transform  
\[  
\hat Wg(t,\xi)=\hat g(\xi)\,t\,\int_{S^2}  
e^{-it\xi\cdot\omega}\,d\omega=\hat  
g(\xi)\,\frac{4\pi\,\sin(|\xi|\,t)}{|\xi|}.   
\]  
Consequently   
\begin{equation}  
I_C^1\leq C\,K^{13/2}\,\nu(K)\,|\log |\delta| |  \,  
\|G^0\|_{L^2}.  
\label{IC1}  
\end{equation}  
Let us now turn to $I_C^2$. We have to define the modified wave operator  
\begin{equation}  
\tilde W_tg(x,v)=t\,\int_{S^2}  
\frac{g(x+t\omega)}{|\frac{v}{\sqrt{1+v^2}}  
\cdot\omega-1|}\,d\omega.\label{modwave}  
\end{equation}  
Note that $\tilde W_t$ enjoys the same regularizing properties as $W$. In particular  
\[  
{\cal F} \nabla_v^k \tilde W_t g (\xi,v)=\hat g(\xi)\,t\,\int_{S^2} \frac{\Phi_k(v,\omega)\,e^{-it\xi\cdot\omega}}{|\frac{v}{\sqrt{1+v^2}}  
\cdot\omega-1|^{k+1}}\,d\omega,  
\]  
for some smooth function $\Phi_k$ of $v$ and $\omega$. Therefore for any $k$  
\begin{equation}  
\|\nabla_x\nabla_v^k \tilde W_t g\|_{L^2_{xv}}\leq C_k\,K^{k+7/2}\,\|g\|_{L^2},\label{tildewl2}   
\end{equation}  
where the $L^2_{xv}$ norm is taken over any regular compact subset of $\R^6$ included in $B(0,K)$. By taking $k$ large enough ($k=2$ for instance) and by Sobolev embedding, one may conclude that  
\begin{equation}  
\|\nabla_x\tilde W_t g\|_{L^2_{x}(L^\infty_v)}\leq C\,K^{15/2}\,\|g\|_{L^2}.\label{tildewlinfty}   
\end{equation}  
Notice now that  
\[\begin{split}  
\left|t\int_{S^2}\frac{G^0(X_t+t\,\omega)  
-G^0(X_t^\delta+t\omega)}{|\dot X_t^\delta\cdot\omega-1|}\,d\omega\right|&=  
|\tilde W_tG^0(X_t,V_t^\delta)-\tilde W_tG^0(X_t^\delta,V_t^\delta)| \\  
\leq |X_t-X_t^\delta|\,(M_t|\nabla_x \tilde W_tG^0(X_t^\delta,V_t^\delta)|&+M_t|\nabla_x  
\tilde W_tG^0(X_t,V_t^\delta)),   
\end{split}\]  
by applying Lemma 3.1 in \cite{Ja}, where we defined the modified  
maximal operator  
\[  
M_sg(x)=\frac{1}{|\delta|+\max_{r\leq s} |X_r-X_r^\delta|}  
\int_{B(x,\max_{r\leq s} |X_r-X_r^\delta|)} \frac{g(z)\,dz}{|z-x|^2},   
\]  
one has  
\[\begin{split}  
I_C^2&\leq C\nu(K)\int_{\Omega_K}\!\!\int_0^T  
\frac{|V_t-V_t^\delta|}{\sqrt{A_\delta(t)}}\left(M_t|\nabla \tilde   
W_tG^0(X_t^\delta,V_t^\delta)|+M_t|\nabla  
\tilde W_tG^0(X_t,V_t^\delta)|\right)\\  
&=I_C^{21}+I_C^{22}.  
\end{split}\]  
The first term can be  easily   bounded  as it is symmetric. By Cauchy-Schwarz, we have   
\[\begin{split}  
I_C^{2,1}\leq C\,\nu(K)\,K^3\,\sqrt{\log  
  \delta}&\,\Bigg(\left(\int_{\Omega_K}\int_0^T  
(M_t|\nabla_x \tilde W_tG^0(X_t^\delta,V_t^\delta)|)^2  
\,dt\,dx\,dv\right)^{1/2}.\\  
\end{split}\]  
Now by Fubini and a change of variable  
\[\begin{split}  
\int_{\Omega_K}\int_0^T  
(&(M_t|\nabla_x \tilde W_tG^0(X_t^\delta,V_t^\delta)|)^2\,dt\,dx\,dv\\  
&\leq \int_0^T  
\int_{B(0,C\,K)}(M_t|\nabla W_tG^0(x,v)|)^2\,dx\,dv\,dt \\  
&\leq C\,\int_0^T  
\int_{B(0,C\,K)}|\nabla_x \tilde W_tG^0(x,v)|^2\,dx\,dv\,dt,  
\end{split}\]  
by the continuity of $M_t$ on $L^2$. By the bound \eqref{tildewl2}, one may bound  
\begin{equation}  
I_C^{2,1}\leq C\,\nu(K)\,K^{11/2}\,\|G^0\|_{L^2}.\label{IC21}  
\end{equation}   
The other term is not symmetric, as it mixes $X_t^\delta$ and $V_t$. It is hence more complicated  
but it can still be handled in a roughly similar way  
\[  
I_C^{2,2}\leq C\,\nu(K)\,K^3\,\sqrt{|\log \delta|}\,  
\left(\int_0^T\!\!\! \int_{\Omega_K}\!\! M_t(|\nabla_x \tilde  
WG^0(X_t,V_t^\delta)|^2)\,dx\,dv\,dt\right)^{1/2}.  
\]  
Now note that  
\[  
M_t g(x)\leq \int_{B(0,K)}  
\frac{|g(z)|\,dz}{(|\delta|+|z-x|)\,|z-x|^2}.
\]  
Hence one gets  
\[\begin{split}  
I_C^{2,2}\leq &C\,\nu(K)\,K^{3}\,|\log \delta|^{1/2}\\  
&\left(\int_0^T  
\int_{(B(0,K))^3} \frac{|\nabla_x \tilde   
WG^0(z,V^\delta_t)|^2}{(|\delta|+|z-X_t|)\,|z-X_t|^2}\,dz\,dx\,dv\,ds\right)^{1/2}.  
\end{split}\]  
Therefore  
\[\begin{split}  
I_C^{2,2}&\leq C\,\nu(K)\,K^{3}\,|\log  
\delta|^{1/2}\\  
&\left(\int_0^T  
\int_{(B(0,K))^3} \frac{\sup_w |\nabla_x \tilde   
WG^0(z,w)|^2}{(|\delta|+|z-X_t|)\,|z-X_t|^2}\,dz\,dx\,dv\,ds\right)^{1/2}.  
\end{split}\]  
Now by the usual change of variables, we get   
\[  
I_C^{2,2}\leq C\,\nu(K)\,K^{9/2}\,|\log  
\delta|\,\|\nabla_x \tilde WG^0\|_{L^2_x(L^\infty_v)}.  
\]  
Using \eqref{tildewlinfty}, one deduces that  
\[  
I_C^{2,2}\leq C\,\nu(K)\,K^{12}\,|\log  
\delta|\,\|G^0\|_{L^2}.  
\]  
Therefore, combining with \eqref{IC1} and \eqref{IC21},   
one finally concludes that  
\begin{equation}  
I_C\leq C\,K^{12}\,\nu(K)\,|\log \delta|\,  
\|G^0\|_{L^2}.\label{IC}  
\end{equation}  
We did not deal with this term in a very subtle manner  
here. However to improve the result, one would need to do $L^1$ or  
at least $L^p$ bounds (instead of $L^2$). 
Note as well that other terms anyway impose the use of the $L^2$ norm.  
\subsubsection{Bound for $I_D$}  
This bound essentially follows the line of $I_A$ in a slightly more complicated  
way.  
  
First of all, one may easily compute the $z$-derivative as  
\begin{equation}  
\nabla_z s_X=\frac{\omega_X}{\omega_X\cdot \dot X_{s_X}-1},\quad \nabla_z  
\omega_X=\frac{1}{s_X}  
\left(\frac{\omega_X\otimes(\dot X_{s_X}-\omega)}{\omega_{X}\cdot  
  \dot X_{s_X}-1}-I\right),\label{diffz}   
\end{equation}  
  and the corresponding formulas for $s_{X^\delta}$ and  
  $\omega_{X^\delta}$.  
Hence  
\[\begin{split}  
\left|\nabla_z\cdot\frac{C\omega_X}  
{s_X\,|\dot X_{s_X}\cdot\omega_X-1|}\right| &\leq  
\frac{C\,|1+\dot V_{s_X}|}{s_X^2\,|\dot X_{s_X}\cdot\omega_X-1|^2}\\  
&\leq \frac{C\,|1+F(s_X,X_{s_X})|}{s_X^2\,|\dot X_{s_X}\cdot\omega_X-1|^2}.  
\end{split}\]  
Inserting the corresponding terms in $I_D$, one finds that  
\[\begin{split}  
I_D\leq &C\,K^2\,\nu(K)\,\int_{\Omega_K}\int_0^T  
\frac{|V_t-V_t^\delta|}{A_\delta(t)}\int_{O_X^t\setminus  
  O^t_{X^\delta}} \frac{|G^0(z)|}{s_X^2}\,|1+F(s_X,  
X_{s_X})|\\  
&+symmetric\ term.   
\end{split}\]  
The only difference with $I_A$ is the additional term  
$F(s_X,X_{s_X})$. Now, note that as before  
\[  
O^t_X\setminus O^t_{X^\delta}\subset\{t-|X_t-X_t^\delta|\leq |z-X_t|\leq t\}.  
\]   
Hence, using spherical coordinates, one gets  
\[\begin{split}  
I_D&\leq C\,K^2\,\nu(K)\,\int_{\Omega_K}\int_0^T  
\frac{|V_t-V_t^\delta|}{A_\delta(t)}\int_{t-|X_t-X_t^\delta|}^t  
|1+F(s,X_s)|\\  
&\qquad\qquad  
\int_{S^2} |G^0(X_t+s\,\omega)|\,d\omega\,ds\,dt\,dx\,dv+symmetric\\  
&\leq C\,K^2\,\nu(K)\,\int_{\Omega_K}\int_0^T  
\frac{|V_t-V_t^\delta|}{A_\delta(t)}M_SG^0(X_t)\,\int_{t-|X_t-X_t^\delta|}^t  
|1+F(s,X_s)|\,ds\\  
&\quad +symmetric,  
\end{split}\]  
with $M_S$ as before the spherical maximal function. By  
Cauchy-Schwartz  
\[\begin{split}  
I_D\leq &C\,K^2\,\nu(K)\Bigg(\int_{\Omega_K}\int_0^T  
\frac{|V_t-V_t^\delta|^2}{A_\delta(t)}\,(M_SG^0(X_t))^2\,dt\,dx\,dv\\  
&\int_{\Omega_K}\int_0^T  
\frac{|X_t-X_t^\delta|}{A_\delta(t)}\int_{t-|X_t-X_t^\delta|}^t  
(1+F(s,X_s))^2\,ds\,dt\,dx\,dv\bigg)^{1/2}.  
\end{split}\]  
The second term is bounded, using Fubini and then changing variables  
to $(x,v)$ from $(X_s,V_s)$, and is less than  
\[  
C\,K^6\,\|F\|_{L^2}^2\leq C\,K^6\,T\,\|F^0\|_{L^2}^2,  
\]   
as the wave operator propagates the $L^2$ norm.  
  
As for the first term, change variables to $(x,v)$ from $(X_t,V_t)$ to  
bound it by  
\[  
\int_{B(0,K)}(M_SG^0(x))^2\,\int_0^T  
\frac{|v-V_t^\delta(V_t^{-1}(x,v))|^2}{\delta^2+\int_0^t  
  |V_s(V_t^{-1}(x,v))-V_s^\delta(V_t^{-1}(x,v))|^2\,ds}\,dt\,dx\,dv  
\]  
Note that by the semi-group property, one still has that  
\[  
|v-V_t^\delta(V_t^{-1}(x,v))|^2=\partial_t \int_0^t  
  |V_s(V_t^{-1}(x,v))-V_s^\delta(V_t^{-1}(x,v))|^2\,ds,  
\]  
and hence the previous term is bounded by  
\[  
\sqrt{-\log |\delta|}\,\int_{B(0,K)}(M_SG^0(x))^2\,dx\,dv\leq C\,K^3\,  
\sqrt{-\log |\delta|}\,\|G^0\|_{L^2},  
\]  
as the spherical maximal function is bounded on $L^2$.   
  
Combining all the estimates, one eventually finds that  
\begin{equation}  
I_D\leq C\,K^{11/2}\,\nu(K)\,\sqrt{|\log \delta|}\,\|G^0\|_{L^2}\,  
\|F^0\|_{L^2}.\label{ID}  
\end{equation}  
Like $I_C$ this term requires the $L^2$ norm of $G^0$. Contrary to  
$I_C$ though, the computation here naturally produces a quadratic term  
in $G^0$ and $F^0$ and one does not see very well how to improve on it.  
\subsubsection{Bound for  $I_E$}  
Applying formula \eqref{diffz}, and recalling that $\phi\leq 1$,  
we can decompose $E$ into  
\[\begin{split}  
|E|\leq & C\,\int_{O^t_X\cap O^t_{X^\delta}}  
|G^0(z)|\,\left|\frac{H(s_X,\omega_X)}{s_X^2\,|\dot  
X_{s_X}\cdot\omega_X-1|^2}-\frac{H(s_{X^\delta},\omega_{X^\delta})}  
{s_{X^\delta}^2\,|\dot  
X_{s_{X^\delta}}\cdot\omega_{X^\delta}-1|^2}\right|\,dz\\  
&+\Bigg|\int_{O^t_X\cap O^t_{X^\delta}} G^0(z)\,\Bigg(  
\frac{\dot V_{s_X}\cdot \omega_{X}}{s_X (\dot  
X_{s_X}\cdot\omega_X-1)^3}\\  
&\qquad\qquad-  
\frac{\dot V_{s_{X^\delta}}\cdot \omega_{X^\delta}}{s_{X^\delta} (\dot  
X_{s_{X^\delta}}  
\cdot\omega_{X^\delta}-1)^3}\Bigg)\,dz\Bigg|=E^1+E^2.  
\end{split}\]  
for some explicit smooth function $H$ whose exact expression is  
unimportant here.  
  
The term $I_E^1$ is bounded by a direct application of Lemma  
\ref{lemIB}  
\begin{equation}  
I_E^1\leq C\,\nu(K)\,K^7\,\sqrt{-\log  
  |\delta|}\,(\|G^0\|_{L^1}+\|G^0\|_{L^p}).\label{IE1}   
\end{equation}  
As for $E^2$, we change back variables to find  
\[\begin{split}  
E^2\leq &\left|\int_0^t (\dot V_s-\dot V_s^\delta) \phi(s)\,ds  
\right|+\Bigg|\int_0^t \dot V_s\,s\, \int_{S^2}  
\Bigg(\frac{\omega\,G^0(X_s+s\omega)}{(1-\dot X_s\cdot\omega)^2}\\  
&\qquad-\frac{\omega\,G^0(X_s^\delta+s\omega)}{(1-\dot  
  X_s^\delta\cdot\omega)^2}  
\Bigg)\,d\omega\Bigg|=E^{2,1}+E^{2,2},   
\end{split}\]  
with  
\[  
\phi(s)=s\,\int_{S^2}G^0(X_s+s\omega)\,\frac{\omega}{(1-\dot  
  X_s\cdot\omega)^2}\,d\omega.   
\]  
The term $E^{2,2}$ is treated in a similar manner as the previous ones  
(note in particular  
that $\dot V_s$ is bounded in $L^2$ by $\|F^0\|_{L^2}$).   
We instead focus on $E^{2,1}$  
and remark that by integration by part  
\[\begin{split}  
I_E^{2,1}\leq &\int_{\Omega_K}\int_0^T \frac{|V_t-V_t^\delta|^2}  
{A_\delta(t)}\,|\phi(t,x,v)|\,dx\,dv\,dt\\  
&+   
 \int_{\Omega_K}\int_0^T \frac{|V_t-V_t^\delta|}  
{A_\delta(t)}\,
\int_0^t |V_s-V_s^\delta|\,|\partial_s   
\phi(s,x,v)|\,ds
\,dt \,dx\,dv.  
\end{split}\]  
Note that  
\[\begin{split}  
\int_{\Omega_K}\int_0^T &\frac{|V_t-V_t^\delta|^2}  
{A_\delta(t)}\,|\phi(t,x,v)|\,dx\,dv\,dt\\  
&\leq C\,K^2  
\int_{\Omega_K}\int_0^T\partial_t \log(  
{A_\delta(t)}|)\,  
t\,\int_{S^2}|G^0(X_t+t\omega)|\,d\omega\\  
&\leq C\,K^2\,|\log |\delta||\,\int_{\Omega_K}\int_0^T\left|  
\partial_t\left(t\,\int_{S^2}|G^0(X_t+t\omega)|  
\,d\omega\right)\right|\,dt\,dx\,dv.  
\end{split}\]  
We remark that  
\[\begin{split}  
\left|\partial_t\left(t\,\int_{S^2} g(X_t+t\omega)  
\,d\omega\right)\right|\leq &C\,K\,\left|t\,\int_{S^2}  
(1,\omega)\,\nabla g(X_t+t\omega)  
\,d\omega\right|\\  
&+\int_{S^2} |g(X_t+t\omega)|  
\,d\omega,  
\end{split}\]  
which implies that this term is bounded in $L^2$ by the $L^2$ norm of  
$g$. Consequently  
\[\begin{split}  
\int_{\Omega_K}\int_0^T &\frac{|V_t-V_t^\delta|^2}  
{A_\delta(t)}\,|\phi(t,x,v)|\,dx\,dv\,dt  
\leq  C\,K^2\,|\log |\delta||\,\|G^0\|_{L^2}.\end{split}  
\]  
As for the other term in $E^{2,1}$, compute  
\[\begin{split}  
|\partial_s\phi|&\leq \int_{S^2}   
\frac{|G^0(X_s+s\omega)|}{(1-\dot X_s\cdot\omega)^2}+s\,  
\left|\int_{S^2}   
\frac{\omega\cdot\nabla G^0(X_s+s\omega)\,\omega}{(1-\dot  
  X_s\cdot\omega)^2}\,ds\right|\\  
&+C\,K^3\,|\dot V_s|\,\int_{S^2} |G^0(X_s+s\omega)|\,ds.  
\end{split}\]  
The first two terms are treated similarly as before. For the last  
term note that by Cauchy-Schwartz  
\[\begin{split}  
&\int_{\Omega_K}\int_0^T \frac{|V_t-V_t^\delta|}{A_\delta(t)}\,\int_0^t  
|V_s-V_s^\delta|\,|\dot V_s|\,\int_{S^2} |G^0(X_s+s\omega)|\,ds\,  
dt\,dx\,dv\\  
&\qquad\leq \left(\int_{\Omega_K}\int_0^T\int_0^t |\dot V_s|^2\,  
\frac{|V_t-V_t^\delta|^2}{A_\delta(t)}\,ds\,dt\,dx\,dv\right)^{1/2}\\  
&\qquad\qquad \left(\int_{\Omega_K}\int_0^T |M_S G^0(X_s)|^2\,  
\frac{\int_0^t |V_s-V_s^\delta|^2\,ds}{A_\delta(t)}\,ds\,dt\,dx\,dv\right)^{1/2},  
\end{split}\]  
with as always $M_S$ the spherical maximal operator. Since $G^0\in  
L^2$, $M_S$ is bounded on $L^2$ and $\dot V_s$ is also bounded in  
$L^2$ by the $L^2$ norm of $G^0$,  
\[\begin{split}  
&\int_{\Omega_K}\int_0^T \frac{|V_t-V_t^\delta|}{A_\delta(t)}\,\int_0^t  
|V_s-V_s^\delta|\,|\dot V_s|\,\int_{S^2} |G^0(X_s+s\omega)|\,ds\,  
dt\,dx\,dv\\  
&\qquad\leq C\, K^3\,\sqrt{-\log|\delta|}\,\|G^0\|_{L^2}^2.  
\end{split}\]  
We conclude that  
\begin{equation}  
I_E\leq C\,\nu(K)\,  
K^7\,\sqrt{-\log|\delta|}\,(\|G^0\|_{L^1}+\|G^0\|_{L^2}+  
\|G^0\|_{L^2}\,\|F^0\|_{L^2}).\label{IE}  
\end{equation}  
\subsubsection{Conclusion of the proof of Lemma \ref{teclem}}   
We combine estimates \eqref{IB0}, \eqref{boundA}, \eqref{IC}, \eqref{ID}  
and \eqref{IE}. Taking the worst terms, one finds the estimate in  
Lemma \ref{teclem}.  
%
\subsection{Conclusion of the proof of Proposition  \ref{propondes}}  
By Lemma \ref{teclem}, one has that  
\[  
\bar Q_\delta(T)\leq C_K\,T\,|\log |\delta||\,(\|G^0\|_{L^1}\,+  
\|G^0\|_{L^2})\,(1+\|F^0\|_{L^2}),  
\]  
for a constant $C_K$ increasing with $K$ (which could be made  
explicit for a given $\nu$).   
  
On the other hand, the definition of $\bar Q_\delta$ yields the very  
obvious bound  
\[  
\bar Q_\delta(T)\leq C_K\,T\,\sqrt{-\log |\delta|}\,\|G^0\|_{W^{1,\infty}}.  
\]  
As $\bar Q_\delta(T)$ depends linearly on $G^0$, one may conclude by  
interpolation that for any fixed $G^0\in L^2$, there exists an  
increasing $\psi_0$  
with  
\[  
\frac{\psi_0(\xi)}{\xi}\longrightarrow 0,\quad\mbox{as}\ \xi\rightarrow  
\infty,  
\]  
such that  
\[  
\bar Q_\delta(T)\leq C_K\,T\,\psi_0(-\log |\delta|).  
\]   
Note that $\psi_0$ depends on $\|G^0\|_{L^2}$ and on the  
equi-integrability of $\hat{G^0}$ in $L^2$. As $F^0\in L^2$, we finally use  
this estimate for $G^0=F^0$ and get  
\[  
Q_\delta(T)\leq \frac{C\,T}{K^2}\,|\log |\delta||+C_K\,T\,\psi_0(-\log  
|\delta|).   
\]   
It only remains to optimize in $K$ by defining  
\[  
\psi(\xi)=\inf_K \left(\frac{C\,\xi}{K^2}+C_K\,\psi_0(\xi)\right).  
\]  
One still has that  
\[  
\frac{\psi(\xi)}{\xi}\longrightarrow 0,\quad\mbox{as}\ \xi\rightarrow  
\infty,  
\]  
and  
\[  
Q_\delta(T)\leq T\,\psi(-\log |\delta|),  
\]  
which concludes the proof of Prop. \ref{propondes}.  
\section{Proof of Theorem  \ref{th1d}}  
\subsection{Well posedness for the corresponding ODE}  
We follow the same steps as for the proof of Th. \ref{ondes}. We study  
the ODE  
\begin{equation}\begin{split}  
&\frac{d}{dt} X(t,x,v)=\alpha(V(t,x,v)),\quad   
\frac{d}{dt} V(t,x,v)=F(t,X(t,x,v),\\  
& X(0,x,v)=x,\qquad V(0,x,v)=v.  
\end{split}  
\label{ode2}  
\end{equation}  
As flows the solutions are again required to satisfy   
\begin{property}  
For any  
  $t\in\R$ the application  
  \begin{equation}  
    (x,v)\in\R\times\R^d\mapsto(X(t,x,v),V(t,x,v))\in\R\times\R^d  
   \label{invertible2}  
  \end{equation}  
  is globally invertible and has Jacobian $1$ at any  
  $(x,v)\in\R\times\R^d$.  
 It also defines a semi-group   
\begin{equation}\begin{aligned}  
&  \forall s,t\in\R,\qquad & X(t+s,x,v)=X(s,X(t,x,v),V(t,x,v)), \\  
 & \mbox{and}\qquad & V(t+s,x,v)=V(s,X(t,x,v),V(t,x,v)).  
\end{aligned}\label{semigroup2}\end{equation}  
\label{prop:Hamilt2}\end{property}  
We  look again  
 at   
\begin{multline*}  
  R_\delta(T)=\log\left(1+\frac{1}{\bar\delta(T,x,v,|\delta|)^2}\left(\sup_{0\leq  
        t\leq T}|X(t,x,v)-X^\delta(t,x,v)|^2  
      \right.\right. \\  
      \left.\left. +\int_0^T|V(t,x,v)-V^\delta(t,x,v)|^2 \:dt\right)\right)  
\end{multline*}  
for any $(x,v)\in \Omega$ a subset of $\R^{d+1}$, and for a function  
$\bar \delta(t,x,v) = \bar \delta(t,x,v,|\delta|^2)   $ to be determined later. $(X,V)$ is a solution to  
\eqref{ode2}  satisfying Property \ref{prop:Hamilt2}  
(or a regularized problem) and $X^\delta,V^\delta$  
verifies  
\begin{equation}  
\begin{split}  
& \mbox{Either}\ (X^\delta,\,V^\delta)\ \mbox{is a solution to  
    \eqref{ode2} (or a regularized  version)}\\ 
&\qquad\qquad\qquad\qquad  \mbox{satisfying Property \ref{prop:Hamilt2}},\\  
& \mbox{Or}\ \exists  
   (\delta_1,\delta_2)\in\R^{1+d}\quad\mbox{with}\ |  (\delta_1,\delta_2) |\leq  
  \delta,\\  
&\qquad\qquad\qquad\qquad  
  (X^\delta,V^\delta)(t,x,v)=(X,V)(t,x+\delta_1,v+\delta_2).    
\end{split}\label{propdelta2}  
\end{equation}  
Theorem  \ref{th1d} is  
a consequence of  
\begin{prop} Assume $F^0\in L^\infty$ and \eqref{mun}.  
For any $\Omega$ compact, any $F^0$  bounded,  
there exists a function $\bar\delta(t,x,v,|\delta|)$, increasing in time with  
$\bar \delta(0,x,v,|\delta|)=|\delta|$, such that  for any 
$(X,V)$ solution to \eqref{ode2} with $F$ given by   
Eq. \eqref{structF}, satisfying Property \ref{prop:Hamilt2},  
and $(X^\delta,V^\delta)$ satisfying \eqref{propdelta2}, for all $T > 0$,   
\[  
R_\delta(T)\leq C\,T\,\|F^0\|_{\infty}\,(\log(1/|\delta|))^{3/4},  
\ \bar \delta(T,x,v,|\delta|)\longrightarrow  
0\ as\ |\delta|\rightarrow 0,\ a.e. (x,v).  
\]\label{prop1d}  
\end{prop}  
  
First of all let us order the $\mu_n$ decreasingly  
\[  
|\mu_0|\geq |\mu_1|\geq \ldots |\mu_n|\ldots  
\]  
Note then that by \eqref{mun}  
\[  
\|F(t,x)\|_\infty\leq \|F_0\|_\infty\,\sum_n |\alpha_n|\leq C\|F^0\|_{\infty}.  
\]  
Therefore defining $\tilde \Omega=\Omega+B(0,T\,\|F_0\|_\infty\,\sum_n  
|\alpha_n|)$, one has that $(X,V)\in\tilde \Omega$ for any  
$(x,v)\in\Omega$ and $t\in [0,\ T]$. The same is of course  
true for $(X^\delta,V^\delta)$.  
  
As before problems occur when the velocity of the particle is close to  
one of the propagation velocities $\xi_n$. So first of all it is  
necessary to control the time that each trajectory spends near one of  
those points.  
  
Denote  
\[\begin{split}  
&\omega(w,\eta,K)=\{(x,v)\in\Omega\ s.t.\ |\{t,\;|\alpha(V(t,x,v))-w|\leq  
\eta\}|\geq K\,\eta\},\\  
&\omega^\delta(w,\eta,K)=\{(x,v)\in\Omega\ s.t.\   
|\{t,\;|\alpha(V^\delta)(t,x,v)-w|\leq   
\eta\}|\geq K\,\eta\}.  
\end{split}\]  
The parameter $\eta$ will be chosen later but will tend to $0$ as  
$|\delta|\rightarrow 0$.  
Then   
\begin{lemma} There exists a constant $C$ depending on  
  $\|F^0\|_{\infty}$ and $\sum_n |\mu_n|$ such that   
\[  
|\omega(w,\eta,K)|\leq \frac{C}{K}.  
\]\label{timecont}  
\end{lemma}   
\noindent{\bf Proof.} Simply write that   
\[   
\int_0^T\int_{\omega(w,\eta,K)} \ind_{ \{ |\alpha(V(t,x,v))-w|\leq  
  \eta  \} }\,dx\,dv\,dt \geq K\,\eta\,|\omega(w,\eta,K)|.  
\]  
On the other hand, using Property \ref{invertible2} and the   
assumption \eqref{alp-prop}   on $\alpha(v)$, we have   
\[\begin{split}  
\int_0^T\int_{\omega(w,\eta,K)} \ind_{ \{  |\alpha(V(t,x,v))-w|\leq  
  \eta  \} }\,dx\,dv\,dt &\leq \int_0^T \int_{\tilde \Omega}  
\ind_{ \{   |\alpha(v)-w|\leq   
  \eta  \}   }\,dx\,dv\,dt\\ 
&\leq  C\,|\tilde\Omega|\,\eta\,T,  
\end{split}\]   
for some constant  $C$.   
 
Finally one  
concludes that   
\[  
|\omega(w,\eta,K)|\leq\,C/K.  
\]  
For $X^\delta,V^\delta$, one uses \eqref{propdelta2} and either they also satisfy \eqref{invertible2} in which case the proof is identical. Or one just has to shift $(x,v)$ by $\delta$ and still follow the same steps.   
\cqfd    
   

Now for $a_n$ (to be fixed later), we  define   
\[  
O_n=\{v,\ |\alpha(v)-\xi_n|\leq a_n\,\eta\}.  
\]  
Then for any $(x,v)$, we decompose the time interval $[0,\ T]$ into  
$I_{x,v}$ and the union $\bigcup_n [t_0^n,\ s_0^n]  
\cup\ldots\cup [t_{k_n}^n,\ s_{k_n}^n]$ with  
\[\begin{split}  
&I_{x,v}=\{t,\ |\alpha(V(t,x,v))-\xi_n|\geq  
a_n\,\eta/2\},\\  
& \sup_{[t_{i}^n,\ s_i^n]} |\alpha(V)-\xi_n|=  
a_n\,\eta,\ \inf_{[t_{i}^n,\ s_i^n]} |\alpha(V)-\xi_n|\leq  
a_n\,\eta/2.  
\end{split}\]  
Similarly one defines $I_{x,v}^\delta$ and the intervals  
$[t_i^{n,\delta},\ s_i^{n,\delta}]$ for $i=0\dots k_n^\delta$.   
Note that $k_n$ and   
$k_n^\delta$ depend  on $(x,v)$.  
  
Define now  
\[  
l_n=\sum_{i=0}^{k_n} (s_i^n-t_i^n),\quad l_n^\delta  
=\sum_{i=0}^{k_n} (s_i^{n,\delta}-t_i^{n,\delta}),  
\]  
and  
\[\begin{split}  
&\omega(\eta,K)=\{(x,v)\in\Omega\ s.t.\ \sum_n\,l_n\geq K\,\eta\},\\  
&\omega^\delta(\eta,K)=\{(x,v)\in\Omega\ s.t.\ \sum_n   
l_n^\delta\geq K\,\eta\}.\end{split}  
\]  
Similarly to Lemma \ref{timecont} one can  deduce  
\begin{lemma} Assume that  $\sum_n a_n<\infty$, then  there exists a constant $C$  
  depending on $\sum_n a_n$ s.t.  
\[   
|\omega(\eta,K)|\leq \frac{C}{K},\quad |\omega^\delta(\eta,K)|\leq \frac{C}{K}.  
\]\label{contln}  
\end{lemma}  
\noindent {\bf Proof.} Simply note that  
\[  
\omega(\eta,K)\subset\big\{(x,v),\quad |\{t,\ V(t,x,v)\in \bigcup_n O_n\}|\geq  
K\,\eta\big\}.   
\]  
Then similarly to the proof of Lemma \ref{timecont}  
\[\begin{split}  
K\,\eta\,|\omega(\eta,K)|&\leq  
\int_{\omega(\eta,K)}\int_0^T \ind_{V(t,x,v)\in\bigcup_n  
  O_n}\,dt\,dx\,dv\\  
&\leq \int_0^T\int_\Omega \ind_{V(t,x,v)\in\bigcup_n  
  O_n}\,dx\,dv\,dt=\int_0^T \int_{\tilde \Omega} \ind_{v\in\bigcup_n  
  O_n}\,dx\,dv\\  
&\leq T\,\sum_n |O_n|\leq T\,\eta\,\sum a_n,    
\end{split}\]  
which shows the result.\cqfd  
  
We are now ready to define $\bar \delta(t,x,v,|\delta|)$. We put  
\begin{equation}\begin{split}  
&\bar \delta(0,x,v,|\delta|)=|\delta|,\quad \partial_t\bar  
\delta(0,x,v,|\delta|) =C\,\|F\|_\infty\,\sum_n (l_n+l_n^\delta).  
\end{split}\label{defbardelta}  
\end{equation}  
Note that from Lemma \ref{contln} and the fact that $\eta$ tends to  
$0$ as $|\delta|\rightarrow 0$, $\bar \delta$ indeed converges to $0$  
for $a.e.\ (x,v)$.   
  
From the computation for $Q_\delta$ in the proof of  
Proposition  \ref{propondes}, one sees that  
\[\begin{split}  
\frac{d}{dt} R_\delta(T)\leq &2\int_0^T  
\frac{V(t)-V^\delta(t)}{A_\delta(t,x,v)}\\  
&\quad \int_0^t (F(s,X(s,x,v))-F(s,X^\delta(s,x,v)))\,ds\,dt\\  
&-2\int_0^T \frac{\partial_t \bar\delta(t,x,v,|\delta|)}{A_\delta(t,x,v)}\;dt,   
\end{split}\]  
with  
\[\begin{split}  
A_\delta(t,x,v)=&\bar \delta^2+\sup_{0\leq s\leq t}  
|X(s,x,v)-X^\delta(s,x,v)|^2\\   
&+\int_0^t |V(s,x,v)-V^\delta(s,x,v)|^2\,ds.  
\end{split}\]  
Compute  
\[\begin{split}  
&\int_0^t (F(s,X(s))-F(s,X^\delta(s)))\,ds\\  
&\qquad=\sum_n \int_0^t  
(F^0(X(s)-\xi_ns)-F^0(X^\delta(s)-\xi_ns))\,\mu_n\,ds\\  
&\ \leq \sum_n \Bigg(C\,\|F\|_\infty\,\mu_n\sum_{i=0}^{k_n} |s_i^n-t_i^n|  
+C\,\|F\|_\infty\,\mu_n  
\sum_{i=0}^{k_n^\delta} |s_i^{n,\delta}-t_i^{n,\delta}|   
\\  
&\qquad\quad+\int_{s\not\in \cup_i(  
  [t_i^n,\ s_i^n]\cup [t_i^{n,\delta},\ s_i^{n,\delta}])}  
(F^0(X(s)-\xi_ns)-F^0(X^\delta(s)-\xi_ns))  
\,\mu_n\,ds\Bigg).\\  
\end{split}  
\]  
  
Now $[0,\ T]\setminus \cup_i  
  [t_i^n,\ s_i^n]\cup [t_i^{n,\delta},\ s_i^{n,\delta}]$ is included  
  in an union  
  of intervals of the form $[s_i^n,\ t_{i+1}^n]$, $[s_i^{n,\delta},  
\ t_{i+1}^{n,\delta}]$, $[s_i^n, t_j^{n\delta}]$ or  $[s_i^{n,\delta},  
\ t_{j}^{n}]$. This depends only on whether $V(s)$ or $V^\delta(s)$ is  
  the first or the last to be such that $|\alpha(V)-\xi_n|=a_n\,\eta$.  
Assume for instance that the corresponding interval is $[s_i^n,\ t_{i+1}^n]$.  
  
Define the transforms  
\[  
\Phi_n(s)=X(s)-\xi_n\,s,\quad \Phi_n^\delta(s)=X^\delta(s)-\xi_n\,s    
\]  
and note that by definition $\Phi_n$ and $\Phi_n^\delta$ are  
invertible on the corresponding interval  $[s_i^n,\ t_{i+1}^n]$. Hence  
denote $S_n(u)$ and $S_n^\delta(u)$ the reciprocal functions.  
  
To bound the next integral, we will make   a change of  variable from $s$ to $u = \Phi_n(s)$ and then go   back  
to the original variable:     
\[\begin{split}  
&\int_{[s_i^n,\ t_{i+1}^n]}  
(F^0(X(s)-\xi_ns)-F^0(X^\delta(s)-\xi_ns))\,ds\\   
&\leq  \int_{\Phi_n(s_i^n)}^{\Phi_n(t_{i+1}^n)}  
F(u)\left(\frac{1}{\alpha(V(S_n(u)))-\xi_n\,S_n(u)}  
-\frac{1}{\alpha(V^\delta(S_n^\delta(u)))-\xi_n\,S_n^\delta(u)}\right)\,du\\  
&\qquad+\|F\|_\infty\;\frac{\sup_s |X(s)-X^\delta(s)|}{a_n\,\eta}\\  
&\leq C\,\frac{\|F\|_\infty}{a_n\,\eta}\;\int_{s_i^n}^{t_{i+1}^n}  
(|V(s)-V^\delta(s)|+|s-S_n^\delta(\Phi_n(s))|)\,ds\\  
&\qquad+\|F\|_\infty\;\frac{\sup_s |X(s)-X^\delta(s)|}{a_n\,\eta},   
\end{split}  
\]  
simply by using the Lipschitz regularity of $\alpha$ in $v$ and of $V,  
V^\delta$ in time. Next notice that  
\[\begin{split}  
|s-S_n^\delta(\Phi_n(s))|&=|S_n^\delta(\Phi_n^\delta(s)-S_n^\delta(\Phi_n(s))|  
\leq \|S_n^\delta\|_{W^{1,\infty}}\;|\Phi_n^\delta(s)-\Phi_n(s)|\\  
&\leq C\,\frac{\sup_s |X(s)-X^\delta(s)|}{a_n\,\eta}.  
\end{split}\]  
And hence deduce that  
\begin{equation}\begin{split}  
&\int_{[s_i^n,\ t_{i+1}^n]}  
(F^0(X(s)-\xi_ns)-F^0(X^\delta(s)-\xi_ns))\,ds\leq  
C\,\frac{\|F\|_\infty}{a_n\,\eta}\\  
&\Bigg(\int_{s_i^n}^{t_{i+1}^n}  
|V(s)-V^\delta(s)|\,ds+(1+|t_{i+1}^n-s_i^n|/(a_n\,\eta))\sup_s  
|X(s)-X^\delta(s)| \Bigg).  
\end{split}\end{equation}  
Introducing this estimate, one obtains  
\[\begin{split}  
\frac{d}{dt} R_\delta(T)\leq & C\,\|F\|_\infty \sum_n \mu_n \int_0^t  
\frac{|V(t)-V^\delta(t)|}{A_\delta(t,x,v)}\;\Bigg(l_n+l_n^\delta\\  
+\frac{1}{a_n\eta}&  
\int_0^t|V(s)-V^\delta(s)|\,ds  
+\frac{k_n+k_n^\delta+t/(a_n\eta)}{a_n\,\eta}\,\sup_{s\leq t}  
|X(s)-X^\delta(s)|\Bigg)   \\  
&-2\int_0^T \frac{\partial_t  
  \bar\delta(t,x,v,|\delta|)}{A_\delta(t,x,v)}\;dt.\end{split}   
\]   
Therefore  
\[\begin{split}  
\frac{d}{dt} R_\delta(T)\leq & C\,\|F\|_\infty\,\sqrt{\log 1/|\delta|}  
\, \sum_n \frac{\mu_n}{a_n\,\eta} (1+k_n+k_n^\delta+1/(a_n\eta))\\  
&+\int_0^T \frac{|V(t)-V^\delta(t)|\sum_n \mu_n  
  (l_n+l_n^\delta)-\partial_t \bar \delta}{A_\delta(t,x,v)}\,dt,   
\end{split}   
\]   
and from the definition \eqref{defbardelta} of $\bar \delta$ and the  
obvious bound $|V-V^\delta|\leq |\delta|+T\|F\|_\infty$, one simply  
gets  
\[  
\frac{d}{dt} R_\delta(T)\leq  C\,\|F\|_\infty\,\sqrt{\log 1/|\delta|}  
\, \sum_n \frac{\mu_n}{a_n\,\eta} (1+k_n+k_n^\delta+1/(a_n\eta)).  
\]  
By the definition of the intervals and the fact that $V$ is Lipschitz  
in time  
\[  
|s_i^n-t_i^n|\geq a_n\,\eta/C.  
\]  
Hence  
\[  
l_n\geq k_n\,a_n\,\eta/C,\ \mbox{or}\ k_n\leq \frac{C\,T}{a_n\,eta}.  
\]  
So finally  
\begin{equation}  
R_\delta(T)\leq  C\,T\,\|F\|_\infty\,\frac{\sqrt{\log 1/|\delta|}}{\eta^2}  
\, \sum_n \frac{\mu_n}{a_n^2}.\label{finalest}  
\end{equation}  
Taking for instance $\eta=(\log 1/|\delta|)^{-1/8}$, one indeed  
concludes the proof of Prop. \ref{prop1d}, provided that it is  
possible to choose the $a_n$ s.t. $\sum_n \mu_n/a_n^2$ is finite.  
  
The only other constraint to satisfy on the $a_n$ is that $\sum_n  
a_n<\infty$. By the bound \eqref{mun}, one may simply choose  
$a_n=n^{-\gamma/2}$. \cqfd   
%

   


\end{document}